\documentclass{llncs}

\usepackage{amssymb}
\usepackage{amsmath}
\usepackage{amsfonts}

\newcommand{\field}[1]{\mathbb{#1}}
\newcommand{\Q}{\field{Q}}
\newcommand{\Z}{\field{Z}}
\newcommand{\F}{\field{F}}
\newcommand{\Proj}{\field{P}}
\newcommand{\Cong}[3]{#1 \equiv #2 \;{\rm mod} \; #3}
\newcommand{\jacobi}[2]{\genfrac{(}{)}{1pt}{}{#1}{#2}}

\newtheorem{conj}{Conjecture}

\begin{document}

\title{Elliptic Curves $x^3 + y^3 = k$ of High Rank}

\titlerunning{Elliptic Curves $x^3 + y^3 = k$ of High Rank}

\author{Noam D. Elkies\inst{1} \and Nicholas F. Rogers\inst{2}}

\authorrunning{Elkies-Rogers}   % abbreviated author list (for running head)

%%%% modified list of authors for the TOC (add the affiliations)
\tocauthor{Noam D. Elkies (Harvard University),
Nicholas F. Rogers (Harvard University)}

\institute{Department of Mathematics, Harvard University, Cambridge, MA 02140\\
\email{elkies@math.harvard.edu}\\
Supported in part by NSF grant DMS-0200687
\and
Department of Mathematics, Harvard University, Cambridge, MA 02140\\
\email{nfrogers@math.harvard.edu}}

\maketitle              % typeset the title of the contribution

\begin{abstract}
We use rational parametrizations of certain cubic surfaces and an
explicit formula for descent via 3-isogeny to construct the first
examples of elliptic curves $E_k: x^3 + y^3 = k$ of ranks $8$, $9$, $10$,
and $11$ over~$\Q$.  As a corollary we produce examples of elliptic
curves over~$\Q$ with a rational $3$-torsion point and rank as high as $11$.
We also discuss the problem of finding the minimal curve $E_k$ of a given rank,
in the sense of both $|k|$ and the conductor of $E_k$,
and we give some new results in this direction.
We include descriptions of the relevant algorithms and heuristics,
as well as numerical data.
\end{abstract}

\section{Introduction}

In the fundamental Diophantine problem of finding rational points on
an elliptic curve $E$, one is naturally led to ask which abelian groups can
occur as the group of rational points $E(\Q)$.  Mordell's theorem
guarantees that $E(\Q)$ is finitely generated, so we have
\[ E(\Q) = E(\Q)_{\mbox{\small tors}} \oplus \Z^r, \]
where $r$ is the rank of $E$.  Mazur's well-known work \cite{mazur}
completely classifies the possibilities for $E(\Q)_{\mbox{\small tors}}$,
but the behavior of the rank remains mysterious.  Part of the ``folklore''
is the conjecture that
there exist elliptic curves with arbitrarily large rank over $\Q$.  But
large rank examples are rare, and the record to date is 24 \cite{mm}.
One might further ask about the distribution of ranks in families of
twists, or with prescribed Galois structure on the torsion subgroup;
there is some evidence to suggest that conditions of this sort do not impose
an upper bound on the rank.

A classical question in number theory is to describe the numbers
$k$ that can be written as the sum of two rational cubes.  This leads
one to study the elliptic curves
\[ E_k: x^3 + y^3 = k \]
for $k \in \Q^{*}$.
Clearly $E_k$ and $E_{k'}$ are isomorphic if $k/k'$ is a cube,
so we can and will restrict our attention to positive cubefree integers $k$.
A Weierstrass equation for $E_k$ is given by $Y^2 = X^3 - 432k^2$, where
\[ X = \frac{12k}{y+x}, \qquad Y = 36k \frac{y-x}{y+x}. \]

As long as $k > 2$, the group $E_k(\Q)_{\mbox{\small tors}}$ is trivial,
so $E_k$ has a nontrivial rational point if and only if its rank
is positive.  The distribution of ranks in this family is not well understood.
Zagier and Kramarz~\cite{zk}
used numerical evidence for $k \leq 70000$ to conjecture that a
positive proportion of the curves $E_k$ have rank at least~$2$; however,
more recent computations by Mark Watkins~\cite{watkins} suggest
that, in fact, a curve $E_k$ has rank 0 or 1 with probability 1.
Still, the following conjecture is widely believed:

\begin{conj}
There exist elliptic curves $E_k$ with arbitrarily large rank over~$\Q$.
\end{conj}

A proof of this conjecture seems beyond the reach of current
techniques.  So for now we content ourselves with constructing
high-rank examples within this family (thereby adding to the
body of supporting evidence), and gathering more data on the
distribution of ranks so as to be able to formulate more precise
conjectures.  The main results of this paper are examples of
curves $E_k$ with rank~$r$ for each $r \leq 11$.  For $r=8,9,10,11$
these are the first curves known of those ranks; for $r=6,7$ our curves
have $k$\/ smaller than previous records, and are proved minimal assuming
some standard conjectures.  For $r \leq 5$
we recover previously known~$k$, and prove unconditionally
that they are minimal.

Throughout, we make use of the fact that the curves $E_k$ are
3-isogenous to the curves
\[ E'_k: uv(u+v) = k \]
or, in Weierstrass form, $E'_k: V^2 = U^3 + 16k^2$, where
\[ U = \frac{4k}{v}, \qquad V = \frac{8ku+4kv}{v}.\]
The isogeny is given by:
\begin{eqnarray*}
\phi: E_k \rightarrow E'_k,
\qquad
(x,y) \mapsto (u,v) = (\frac{y^2}{x}, -\frac{k}{xy}).
\end{eqnarray*}
The dual isogeny, with respect to the Weierstrass equations
for $E_k$ and $E'_k$, is
\[ \hat{\phi}: E'_k \rightarrow E_k,
\qquad
(U,V) \mapsto (X,Y) =
\left(\frac{U^3+64k^2}{U^2}, \frac{V(U^3-128k^2)}{U^3}\right).
\]

Applying Tate's Algorithm \cite{tate} to the curves $E'_k$,
we find that a minimal Weierstrass form for $E'_k$ is given by
\[ Z^2 = W^3 + \frac{k^2}{4} \qquad (W,Z) = (\frac{U}{4}, \frac{V}{8}) \]
in the case that $k$ is even, and 
\[
Z^2 + Z = W^3 + \frac{k^2 - 1}{4} \qquad
(W,Z) = (\frac{U}{4}, \frac{V-4}{8})
\]
in the case that $k$ is odd.  The primes of bad reduction for $E'_k$
are the primes dividing $k$ and the prime 3.
For a prime factor~$p$ of~$k$, the Kodaira type at~$p$ is
$\mbox{IV}^*$ if $p^2 \mid k$ and $\mbox{IV}$ if $p \| k$.
If $3 \nmid k$, the Kodaira type at $3$ is $\mbox{III}$
if $k \equiv \pm 2 \mod{9}$, and $\mbox{II}$ otherwise.
It follows that the conductor of $E'_k$ is given by the formula
\[ N(E'_k) = \prod_{p \mid 3k} p^{2 + \beta_p} \]
where $\beta_p = 0$ if $p \neq 3$, and $\beta_3 = 0$, $1$, or $3$
for $k \equiv \pm 2 \mod{9}$, $k \equiv \pm 1, \pm 4 \mod{9}$, or $3 \mid k$
respectively.

The curves $E'_k$ have the rational 3-torsion points $(U, V) = (0, \pm 4k)$.
Since the rank is an isogeny invariant, we produce as a corollary to our work
examples of elliptic curves $E'_k$ with a rational \hbox{$3$-torsion} point
and rank as high as~$11$.
Curiously, there are no other known elliptic curves over~$\Q$ with a rational
\hbox{$3$-torsion} point and rank greater than~$8$~\cite{dujella}.

In section 2 we describe the geometric underpinnings of our search
technique, which heavily uses rational parametrizations of various
cubic surfaces, the points of which correspond to pairs of (usually
independent) points on the curves $E_k$.
Section 3 gives a formula for an upper bound on the rank of $E_k$,
using descent via 3-isogeny.  Section 4 describes some of
the specific algorithms we used to produce examples of $E_k$ with high
rank.  Finally, we give our numerical results in Section 5.

\section{Cubic Surfaces}

The most na\"{\i}ve approach to constructing curves $E_k$ of high rank is
to enumerate small points on the curves $E_k$, which can be accomplished
via the simple observation that a point on some curve $E_k$
corresponds to a pair of whole numbers $(x, y)$ so that the cubefree part
of $x^3 + y^3$ (that is, the unique cubefree integer $s$ such that
$(x^3 + y^3)/s$ is a perfect cube) is $k$.
The second author used essentially this
approach to find the first known $E_k$ of rank~$7$.  By incorporating
some more sophisticated techniques, such as 3-descent (see below),
this approach could yield curves with rank as high as~$8$.
The weakness of this method is that the number of such points
up to height $H$ grows as $H^2$,
most of which lie on curves of rank 1 and waste our time and/or memory.

We can reduce this $H^2$ to $H^{1+\epsilon}$ by considering only
curves $E_k$ together with a pair of points,
which correspond to points on the cubic surface
\[ S_1: w^3 + x^3 = y^3 + z^3 \]
other than the trivial points on the lines
$w+x=y+z=0$, $w+y=x+z=0$, $w+z=x+y=0$.
(This pairs-of-points idea is also used in~\cite{elkies-watkins}
to produce elliptic curves with high rank and smallest conductor known.)
The cubic surface whose points correspond to pairs of points
on the isogenous curves $E'_k$,
\[ S_2: wx(w+x) = yz(y+z),\]
and the ``mixed'' cubic surface
\[ S_3: wx(w+x) = y^3 + z^3,\]
are also fruitful.
Each of these cubic surfaces is smooth, and thus rational
in the sense that it is birational to $\Proj^2$ over $\overline{\Q}$;
in fact, each has a rational parametrization defined over $\Q$.

A parametrization of $S_1$ was found by the first author~\cite{cubic_param},
and a parametrization of $S_2$ follows fairly quickly:
there is an obvious isomorphism between $S_1$ and $S_2$,
defined {\it a priori} over $\Q(\sqrt{-3})$, but which actually
descends to $\Q$.  To parametrize $S_3$, we used the following
more general approach, provided by Izzet Coskun~\cite{izzet}.

Let $S$ be a cubic surface defined over $\Q$, and suppose
$L_1$ and $L_2$ are disjoint lines on $S$, with $L_3$ a third line meeting
both.  Then there is a 3-dimensional space of quadrics that vanish on
this set of lines.  Use a basis of this space to map $S$ into $\Proj^2$;
the inverse map will be a parametrization of $S$.  The parametrization
so obtained is defined over $\Q$ if all of the $L_i$ are rational, or
if $L_3$ is rational and $L_1, L_2$ are Galois conjugate.  This construction
realizes $S$ as $\Proj^2$ blown up at six points; the six blown-down
lines are the six lines (other than $L_i$) that meet exactly two
of the three lines $L_1$, $L_2$, $L_3$.

On both $S_1$ and $S_3$, the relative paucity of lines defined over
$\Q$ means that, up to automorphism of the surface, there is only one
choice for the configuration $L_1$, $L_2$, $L_3$.  Thus for $S_1$ the
parametrization obtained with this technique,
\begin{eqnarray*}
(w:x:y:z) &=& (t^3-2t^2s-2tsr+ts^2+r^2s-r^3-rs^2\\
		&&:-t^3-2t^2s+ts^2-2tsr+2rs^2-s^3+r^3-2r^2s\\
		&&:2t^3+3t^2r-2t^2s-2tsr+3tr^2+ts^2+2rs^2-s^3+r^3-2r^2s\\
		&&:-2t^3+t^2s-3t^2r-3tr^2+4tsr-2ts^2-r^3+r^2s-rs^2),
\end{eqnarray*}
is equivalent to the one in~\cite{cubic_param}, in the sense that
one can be obtained from the other by composing an automorphism of $S_1$
with a projective linear transformation of $\Proj^2$.
For $S_3$ we obtain the parametrization
\[ (w:x:y:z) = (r^3 - s^3:s^3 + t^3:r^2s - s^2t + t^2r:r^2t - s^2r - t^2s). \]
In order to compute the $k$ for which a particular point on $S_3$
corresponds to a pair of points on $E'_k$ and $E_k$, we must
find the cubefree part of $wx(w+x)$, which we do by factoring this
number.  It is therefore useful that $wx(w+x)$, which is a polynomial
of degree 9 in $r$, $s$, and $t$, decomposes as a product
of three linear and three quadratic factors.  By contrast,
the factorization of $w^3 + x^3$ in the parametrization of $S_1$ is
as a product of one linear, two quadratic, and one quartic factor;
the difficulty of factoring the value of this quartic at $(r,s,t)$
severely limits the usefulness of the $S_1$ parametrization.

On $S_2$, there are, up to automorphism, seven different configurations of
lines $L_1$, $L_2$, $L_3$.  Four of them lead to parametrizations
where $wx(w+x)$ has five linear and two quadratic factors;
the parametrization of $S_2$ mentioned above is of this type.
The other three lead to factorizations into three linear and three
quadratic polynomials.  There is not a significant computational
advantage for one factorization over another, but we do mention here
a rather elegant parametrization of $S_2$ obtained in this way:
\[ (w:x:y:z) = (-r^2s + s^2t:r^2s - rt^2:-r^2t + st^2:rs^2 - st^2).\]

\section{Descent via 3-isogeny}

A powerful tool for obtaining upper bounds for the ranks of the curves
$E_k$ is descent, since these curves admit
the aforementioned \hbox{$3$-isogeny} with $E'_k$.  An analysis of the
descent for these curves first appeared in \cite{selmer}.  What
follows is essentially a simplification of the formula given there.

Let $k = \prod p_j^{\epsilon_j}$, where $\epsilon_j = 1$ or 2, and let
$q_i$ be the primes dividing $k$ with $\Cong{q_i}{1}{3}$.  Now define
a matrix $A$ over $\F_3$ by:
\[ \jacobi{p_j^{\epsilon_j}}{q_i}_3 = \rho^{A_{ij}} \]
if $p_j \neq q_i$, and $A_{ij} = -\sum_{\ell : p_\ell \neq q_i} A_{i\ell}$ if
$p_j = q_i$ (equivalently, the rows of $A$ sum to zero).  Here
$\jacobi{\cdot}{q}_3$ denotes the cubic residue symbol mod~$q$;
note that for each $\Cong{q}{1}{3}$, there are two choices
for this cubic residue symbol, but they lead to proportional rows $A_i$.
If $\Cong{k}{\pm 1 {\rm \ or\ } 0}{9}$,
add an additional row corresponding to cubic characters mod 9
for $p_j$; if $9 | k$, the entry corresponding to $p_j = q_i = 3$
is defined as in general when $p_j = q_i$.

If we let $\phi: E_k \rightarrow E'_k$ and $\phi': E'_k \rightarrow
E_k$ denote the relevant $3$-isogenies, then the row and column null
spaces of $A$ correspond, respectively, to the $\phi$- and
$\phi'$-Selmer groups of $E_k/\Q$ and $E'_k/\Q$.  We can conclude
that, after taking the 3-torsion of $E'_k$ into account,
\[ \text{rank}(E_k) \leq \#\text{rows} + \#\text{columns} - 2 \cdot
\text{rank}(A) - 1. \]

\section{Computational Techniques}

An application of the explicit formula for descent via 3-isogeny is another
technique for searching for curves of large rank, which tends to be
more effective in finding the minimal $k$ such that $E_k$ has a given
rank (it was actually this
technique that produced the current rank 9 record).  The idea is to
enumerate all possible $k$ less than some given upper bound which have
a sufficiently high 3-Selmer bound, and then to search for points on
these curves.  To do this, we recursively build up products of primes,
and at each stage compute the portion of the matrix $A$ corresponding
to the primes chosen so far.  Of course, the diagonal entries remain
in some doubt, since they depend on the whole row of $A$; still, at
each stage a lower bound for the rank of the matrix $A$ can be
computed, and used to give a lower bound on the number of primes still
needed.  In this way one can vastly reduce the search space.
A similar approach can be used to enumerate candidate curves
whose conductor is smaller than a given bound,
with the formula for the conductor given in the Introduction.

It is also important to have a way to guess which curves are the most
promising before committing to a lengthy point search.  The most
important tool in this regard is provided by a heuristic argument
suggesting that high rank curves should have many points on their
reductions modulo~$p$, or more specifically,
\[ \prod_{p \leq x} \frac{\#E(\F_p)}{p+1} \sim (\log x)^r \]
where $r$ is the rank of $E$.  This formula was conjectured by Birch
and Swinnerton-Dyer in \cite{bsd2}, and the idea of using it to find
elliptic curves of high rank is due to Mestre \cite{mestre}.

Note that for the curves $E_k$ it is only
useful to consider primes $\Cong{p}{1}{3}$, since $E_k$ is
supersingular mod~$p$ when $\Cong{p}{2}{3}$.  Furthermore, computing
$\#E_k(\F_p)$ is quite fast:  modulo each prime $\Cong{p}{1}{3}$,
there are only three isomorphism classes of $E_k$,
corresponding to the three cubic residue symbols mod~$p$.
We compute the $a_p$'s for each of these isomorphism classes
once for all, and then to find the $a_p$ for a given curve $E_k$,
we need only compute the cubic residue symbol of $k$ mod~$p$.

In the end, though, we must still search for points on
curves we suspect of having high rank.  But here, too, there are
improvements over the most na\"{\i}ve approach.  As noted above, points on
$E_k$ correspond to pairs of whole numbers $(x,y)$ such that
$k = d^{-3}(x^3 + y^3)$.  Of course, we may further assume that $x$ and
$y$ are relatively prime.  It follows that $\gcd(x + y, x^2 - xy + y^2)$
is either~$1$ or~$3$, whence $x + y$ must be a factor of~$3k$
times a perfect cube.  For each $x + y$ and $d$, we must simply decide
if the resulting quadratic equation has a rational solution.
Furthermore, we can use local conditions
to reduce the number of possible $x + y$ we must consider
(this is closely related to the descent described in the last section).
A similar approach works for a point search on $xy(x+y) = k$.

\section{Results}

Here we list the minimal known $k$ such that $E_k$ has rank $r$
for each rank $r\leq 11$, as well as the minimal known conductor
of a curve $E_k$ of rank $r$ with $r \leq 8$.  We include notes where
relevant, and $r$ independent points for each of the new record curves.
The points are listed on the minimal Weierstrass equation for $E'_k$,
as described in the Introduction.  To transfer points back to the curve $E_k$
one may use the dual isogeny $\hat{\phi}$ described there.

The following records for ranks up to 5 are known to be minimal
(the proof for ranks 4 and 5 seems to be new); for ranks 6 and~7,
they are minimal provided the weak Birch and Swinnerton-Dyer Conjecture
and the Generalized Riemann Hypothesis are true for all $L(E_{k'},s)$
with $k' < k$.  The records for
ranks 8 through 10 are likely to be minimal.  In each case, we use the
approach described in the last section to enumerate all of
the smaller $k$ with a sufficiently large 3-Selmer group.  For each of
them we compute a partial product of $L(E_k,1)$ over
the first 1000 or so primes; a large partial product should correspond
to high rank.  In each case of rank 8 through 10, the record $k$
significantly distinguished itself from all smaller $k$.
Note that for large $r$ our record value of $k$ tends to have
considerably more prime factors congruent to~1 mod~3 than to~2 mod~3.

\vspace{0.5cm}

\begin{tabular}{c|l}
rank & \;$k$\;\mbox{(minimal known with $E_k$ of given rank)} \\
\hline
0 & \;$1$ \\
1 & \;$6 = 2 \cdot 3$ \\
2 & \;$19$ (prime) \\
3 & \;$657 = 3 \cdot 3 \cdot 73$ \\
4 & \;$21691 = 109 \cdot 199$ \\
5 & \;$489489 = 3 \cdot 7 \cdot 11 \cdot 13 \cdot 163$ \\
6 & \;$9902523 = 3 \cdot 73 \cdot 103 \cdot 439$ \\
7 & \;$1144421889 = 3 \cdot 13 \cdot 19 \cdot 41 \cdot 139 \cdot 271$ \\
8 & \;$1683200989470 = 2 \cdot 3 \cdot 5 \cdot 7 \cdot 11 \cdot 13 \cdot
17 \cdot 29 \cdot 41 \cdot 47 \cdot 59$ \\
9 & \;$349043376293530 = 2 \cdot 5 \cdot 37 \cdot 41 \cdot 53 \cdot 73
\cdot 1231 \cdot 4831$ \\
10 & \;$137006962414679910 = 2 \cdot 3 \cdot 5 \cdot 7 \cdot 23 \cdot 31
\cdot 37 \cdot 43 \cdot 83 \cdot 109 \cdot 151 \cdot 421$ \\
11 & \;$13293998056584952174157235$ \\
   & \;$= 3 \cdot 5 \cdot 7 \cdot 13 \cdot 19
\cdot 23 \cdot 31 \cdot 43 \cdot 59 \cdot 61 \cdot 73 \cdot 79 \cdot
103 \cdot 109 \cdot 157 \cdot 457$
\end{tabular}

\vspace{0.5cm}

The $k$ in the following chart are known to correspond to curves $E_k$
of minimal conductor for $k \leq 3$.  For ranks $4$, $5$, and $6$,
they are minimal provided the weak Birch and Swinnerton-Dyer Conjecture
and the Generalized Riemann Hypothesis are true
for all $L(E_{k'},s)$ with $k' < k$.
The records for ranks $7$ and $8$ are likely to be minimal;
as above, Mestre's heuristic distinguishes them markedly from all curves
of smaller conductor.  It is striking that for all ranks less than $8$,
the $k$ corresponding to minimal conductor are squarefree away from $3$.
Since divisibility by $p^2$ as opposed to $p$ in $k$ does not alter
the value of the conductor of $E_k$, one might expect that for high rank,
the $k$ corresponding to the curve of minimal conductor would have many
square factors.

\vspace{0.5cm}

\begin{tabular}{c|l}
rank & \;$k$\;\mbox{($E_k$ of given rank and minimal known conductor)} \\
\hline
0 & \;$1$ \\
1 & \;$9 = 3^2$ \\
2 & \;$19$ (prime) \\
3 & \;$657 = 3 \cdot 3 \cdot 73$ \\
4 & \;$34706 = 2 \cdot 7 \cdot 37 \cdot 67$ \\
5 & \;$763002 = 2 \cdot 3^2 \cdot 19 \cdot 23 \cdot 97$ \\
6 & \;$24565833 = 3^2 \cdot 17 \cdot 307 \cdot 523$ \\
7 & \;$1144421889 = 3 \cdot 13 \cdot 19 \cdot 41 \cdot 139 \cdot 271$ \\
8 & \;$23381862574950 = 2 \cdot 3^2 \cdot 5^2 \cdot 11 \cdot 19 \cdot 23 \cdot
31 \cdot 83 \cdot 4201$ \\
\end{tabular}

\vspace{0.5cm}

We conclude by listing independent points on the record curves $E'_k$
for ranks 6 through 11, all of which were newly found using methods
described here.  A program that implements LLL reduction on the lattice
of points of $E_k$, provided by Randall Rathbun \cite{rathbun},
was used to reduce their heights where possible.

\samepage{
\[ k = 9902523,\; r = 6 \]
6 independent points on the Weierstrass minimal curve for $E'_k$:
\[ y^2 + y = x^3 + 24514990441382 \]
($100092$, $32051170$), ($-6798$, $4919434$), ($-22338$, $3656314$),\\
($43672$, $10383069$), ($-11988$, $4774114$), ($126720$, $45380386$).
}

\samepage{
\[ k = 24565833,\; r = 6 \]
6 independent points on the Weierstrass minimal curve for $E'_k$:
\[ y^2 + y = x^3 + 150870037745972 \]
($-37656$, $9872932$), ($86292$, $28167835$), ($187270$, $81966083$),\\
($-32058$, $10859260$), ($-39798$, $9372019$), ($236572$, $115719195$).
}

\samepage{
\[ k = 1144421889,\; r = 7 \]
7 independent points on the Weierstrass minimal curve for $E'_k$:
\[ y^2 + y = x^3 + 327425365005582080 \]
($267748$, $588744383$), ($1235988$, $1488490048$), ($-333330$, $538877944$),\\
($-648774$, $233133847$), ($-422760$, $501863680$), ($5104008$, $11545190143$),\\
($-688974$, $19483912$).
}

\samepage{
\[ k = 1683200989470,\; r = 8 \]
8 independent points on the Weierstrass minimal curve for $E'_k$:
\[ y^2 = x^3 + 708291392738196762720225 \]
($-88860785$, $81396479060$), ($-87348261$, $204569627262$),\\
($-63256830$, $674665720365$), ($-40588401$, $800890328532$),\\
($101707060$, $1326794024465$), ($-35705670$, $814107128835$),\\
($-44793980$, $786391914385$), ($8308684440$, $757353550270065$).
}

\vspace*{2ex}

In fact, the curve $E'_k$ for $k=1683200989470$ has the remarkable
property that the Diophantine equation $xy(x+y) = k$ actually has 8
{\it integral} solutions, namely
$(11,391170)$, $(533,55930)$, $(770,46371)$, $(1003,40467)$,
$(2639,23970)$, $(6970,12441)$, $(7293,1197)$, $(8555, 10387)$.
These 8 solutions, considered as rational points on~$E'_k$,
are independent.

\samepage{
\[ k = 23381862574950,\; r = 8 \]
8 independent points on the Weierstrass minimal curve for $E'_k$:
\[ y^2 = x^3 + 136677874368461861091875625 \]
($1826174700$, $78910161016275$), ($794409100$, $25259021976275$),\\
($-259483950$, $10918164759975$), ($499986216$, $16176140969139$),\\
($-503804925$, $2966885463150$), ($798185799$, $25400836633032$),\\
($165873591$, $11884516327764$), ($215137494$, $12109307517153$).
}

\samepage{
\[ k = 349043376293530,\; r = 9 \]
9 independent points on the Weierstrass minimal curve for $E'_k$:
\[ y^2 = x^3 + 30457819633596695100179965225: \]
($-734843410$, $173381106196815$), ($5130038900$, $406775844709485$),\\
($-2676929565$, $106184137901590$), ($690947990$, $175464197227935$),\\
($291207945620$, $157146639625792365$), ($25120488440$, $3985280944128435$),\\
($-872639080$, $172607374924365$), ($2918890200$, $235215923409485$),\\
($-102315705$, $174518619468760$).
}

\samepage{
\[ k = 137006962414679910,\; r = 10 \]
10 independent points on the Weierstrass minimal curve for $E'_k$:
\[ y^2 = x^3 + 4692726937524378378756566939402025: \]
($-135797482140$, $46781315964225555$), ($-150436201545$, $35891470127810220$),\\
($-42200591214$, $67952721291406041$), ($2327642247924$, $3551854243978575507$),\\
($5504535148140$, $12914782107290941395$), ($140506152430$, $86409469562070095$),\\
($397507563420$, $259814927561209005$), ($7162660587075$, $19169656506442936830$),\\
($73148794740$, $71303068454026605$), ($-102758626586$, $60063833881519937$).
}

\samepage{
\[ k = 13293998056584952174157235,\; r = 11 \]
11 independent points on the Weierstrass minimal curve for $E'_k$:
\[ y^2 + y = x^3 + 44182596082121121317135170025680399046545625711306:\]
($-30156002278649820$, $4093799681127459731025817$),\\
($11364087102067560$, $6756491872572362690626342$),\\
($-20835788771691894$, $5927660006237675713476241$),\\
($1134264920569989390$, $1208031685828825118221478017$),\\
($8907565209691176834$, $26585114133655761890666064910$),\\
($111849199886121334$, $37992674604901443769570910$),\\
($11724873521668020$, $6767159346634715672034457$),\\
($-138658831412368575/4$, $12719819443574268333325811/8$),\\
($165971060901522240$, $67941788876402816577138982$),\\
($994768217796990$, $6647073075327662243966017$),\\
($532896351059436225/16$, $576457310785324883248677823/64$).
}

\end{document}